\documentclass[12pt,a4paper]{article}
\usepackage{a4}
\usepackage[dvips]{graphicx}
\usepackage{amsfonts,amssymb,amsmath}
\begin{document}

\newtheorem{theorem}{Theorem}[section]
\newtheorem{thrm}{Theorem}[section]
\newtheorem{remark}{Remark}[section]
\newcommand{\qed}{\begin{flushright}$\Box$\end{flushright}}
\def\RR{\mathbb{R}}
\renewcommand{\theequation}{\thesection.\arabic{equation}}
\newtheorem{rmrk}{Remark}[section]

\begin{center}
\Large \bf Explicit solutions for a nonlinear model of financial derivatives
\end{center}
\begin{center} {L. A. Bordag${}^*$\footnote{e-mail: Ljudmila.Bordag@ide.hh.se}\\[5pt]
{ \it Halmstad University, Box 823, 301 18 Halmstad, Sweden}\\[5pt]
A. Y. Chmakova ${}^*$}\footnote{e-mail: chmakova@math.tu-cottbus.de}\\[5pt]
 {\it Fakult{\"a}t Mathematik, Naturwissenschaften und Informatik\\
  Brandenburgische Technische Universit{\"a}t Cottbus\\
  Universit{\"a}tsplatz 3/4, 03044 Cottbus, Germany}\\[5pt]
\end{center}

\begin{abstract}
  Families of explicit solutions are found to a nonlinear Black-Scholes
  equation which incorporates the feedback-effect of a large trader in case of
   market illiquidity.  The typical solution of these families will have a payoff
  which approximates a strangle. These solutions were used to test numerical schemes for solving a nonlinear
  Black-Scholes equation.  

\end{abstract}
\noindent
{\bf Key words and phrases:} Black~-~Scholes model, illiquidity,
nonlinearity,\\ explicit solutions\\
{\bf AMS classification:}
35K55, 22E60, 34A05 \\ [5pt]

\setlength{\parindent}{0cm}

\newpage

\section{Introduction}

Standard option pricing theory uses a number of basic
assumptions including the assumptions of symmetric information, of
complete and frictionless markets, as well as the assumption 
that all participants act
as price takers. Recently a series of papers appeared in which one or
more of these assumptions have been relaxed; \cite{WilmottWhalleyj}, 
\cite{WhalleyWilmottg}, \cite{DuffieFlemingSonerZariphopoulouab},
\cite{WhalleyWilmottq}, \cite{Freyj} and
\cite{FouquePapanicolaouSircarSolnac}
are representative examples of this work.
The turbulence on financial
markets such as the events surrounding the collapse of LTCM in 1998
 have made market liquidity an issue of high concern for
investors and risk managers and have triggered a lot of academic research;
see for instance  \cite{SchonbucherWilmotti}, \cite{Frey:perfect}, 
\cite{SchonbucherWilmotta}. 
In illiquid markets an attempt to buy/sell a large amount of an asset will
affect its price so that the assumption that investors act act as price takers 
cannot be maintained.

The purpose of this paper is to
investigate the evaluation of an option hedge-cost under relaxation of the
price-taking assumption. For our analysis we use the framework proposed by
Frey in \cite{Frey:perfect}, \cite{Frey:illicuidity}. He 
developed a model of market illiquidity 
 describing the asset price dynamics which result if a large 
trader chooses a
given stock-trading strategy $(\alpha_t)_t$.
The resulting stock-price dynamics have the following natural property: if the
large trader buys (sells) stock, i.e., if 
 ${\rm d} \alpha_t >0$ (${\rm d}\alpha_t <0$) 
the stock price rises (falls). If the position of the large trader is
unchanged,
the stock price $S_{t}$ follows  standard geometric Brownian motion with
constant volatility $\sigma$. Formally, Frey models stock price dynamics by
the following 
 stochastic differential equation
\begin{equation}
{\rm d} S_t = \sigma S_{t-}{\rm d} W_t + \rho S_{t-}{\rm d} \alpha_t, 
\label{stomo}
\end{equation}
where $ W_t$
 is a standard Brownian motion and
$S_{t-}$ denotes the left limit $\lim _{s \to t, s<t} S_t$. 
In (\ref{stomo})
$\rho$ is the market illiquidity parameter with $0\le \rho $. 
The value $1/(\rho S_t)$ is called depth of the market at time $t$.
Note that in the model (\ref{stomo}) the parameter $\rho $ is a characteristic
of the market and
does not depend on the payoff of the hedged derivatives.
If $\rho \to 0$ then  (\ref{stomo}) reduces to the Black--Scholes model.
We concentrate our investigations on the nontrivial case $\rho \ne 0$.

Consider the problem of hedging a terminal-value claim with maturity $T$ and
payoff $h(S)$ in the model (\ref{stomo}).
 As shown in \cite{Frey:perfect}, \cite{Frey:illicuidity},
 the feedback-effect
leads to a nonlinear version of the Black--Scholes partial differential
equation for a hedge cost $u(S,t)$ of the claim,
\begin{eqnarray} \label{urav}
u_t+\frac{\sigma^2 S^2}2\frac{u_{SS}}{(1-\rho S u_{SS})^2}=0,
\end{eqnarray}
with terminal condition $u(S,T)=h(S)$. The variable $S$ denotes the price 
of the underlying asset and $t$ is the time variable. The equation above 
is studied for the variables $S$ and $T$ in the intervals
\begin{equation} \label{inter}
S \ge 0, ~~ t \in [0, T],~~ T >0.
\end{equation}
Similar equations in related models were obtained by a number of authors 
see for instance \cite{Frey:market}, \cite{Frey:perfect},
\cite{SircarPapanicolaou}, \cite{SchonbucherWilmotti}, 
\cite{SchonbucherWilmotta}, \cite{WhalleyWilmottg}, \cite{WhalleyWilmottq}.

Frey and co-authors , \cite{Frey:illicuidity}, \cite{FreyPatie} studied
equation (\ref{urav}) under constraint and did some numerical
simulations.  Our goal is to investigate this equation using
analytical methods.

We study the model equation (\ref{urav}) using methods of Lie
group theory in Section \ref{liesec2}.
 Using the symmetry group we reduce the
partial differential equation (\ref{urav}) to an ordinary
differential equation in Section \ref{scava}. We obtain
 nontrivial explicit solutions for this case. We prove that
the explicit solutions approximates strangles with corresponding
payoffs (see Section \ref{fainv}). Further,  in Section \ref{prop5} we 
study different properties
of the  obtained solutions. The existence of
nontrivial explicit solutions allows us to test different
numerical methods usually used to calculate hedge-costs of derivatives.
The best results are achieved by the completely implicit method.
The validated numerical
 scheme was used to calculate option hedge-costs in case of calls and
bull-price-spreads.

\section{Lie group symmetries}\label{liesec2}

In this section we study the symmetry properties of equation
(\ref{urav}) and obtain the complete description of the corresponding
Lie algebra, the associated Lie group and a list of functionally
independent invariants.

Let us study the nonlinear part of this equation.
The denominator in the second term of this equation will
be equal to zero if the function $u(S,t)$ satisfies the equation
\begin{eqnarray} \label{denom}
1-\rho S u_{SS}=0.
\end{eqnarray}
The solution of this equation is a function $u_0(S,t)$,
\begin{eqnarray} \label{nuli}
u_0(S,t)=\frac{1}{\rho} S \ln S + S c_1(t) + c_2 (t),~~\rho \ne  0,
\end{eqnarray}
where the functions $c_1(t)$ and $c_2(t)$ are arbitrary functions of the
variable $t.$ From now on we assume that the denominator in the
second term of equation (\ref{urav}) is  not
identically zero, i.e., the function $u(S,t)$ is not equal to the function
$u_0(S,t)$ (\ref{nuli}) except in a discrete set of points.

\vspace{2pt}
We introduce the necessary notations connected with the Lie group theory.
Besides the classical work \cite{Lie} our notations follow \cite{Gaeta} and,
especially with respect to the invariants, to Ovsiannikov \cite{Ovsiannikov} and Olver \cite{Olver}.
We
introduce the two-dimensional space $X$ of independent
variables $(S,t) \in X$ and a one-dimensional space of the dependent
variables $u \in U.$ Then we consider the space $U_{(1)}$ of the first
derivatives of the variable $u$ on $S$ and $t$,
i.e., $(u_S,u_t) \in U_{(1)}.$ Analogously we introduce  the space $U_{(2)}$ of
the second order derivatives $(u_{SS}, u_{St}, u_{tt}) \in U_{(2)}.$
Let $M=X \times U$ be the Cartesian product of pairs $(x,u)$ with $x=(S,t) \in X,~~u\in U$.

The second order jet bundle $ M^{(2)}$ of the base space $M$
has the form
\begin{equation}
M^{(2)} = X \times U \times U_{(1)} \times U_{(2)} \label{jet}.
\end{equation}
We label the coordinates in the space 
 $M^{(2)}$ by $w=(S,t,u,u_S,u_t,u_{SS}, u_{St},
u_{tt}) \in M^{(2)}.$
The second order jet bundle $ M^{(2)}$ has a natural contact structure
 (see \cite{Ovsiannikov}, \cite{Olver}, \cite{Stephani}, \cite{Gaeta}, 
\cite{Ibragimov}).
Our differential equation (\ref{urav}) is of order two and 
in the context of
the second order jet bundle $M^{(2)}$ it should be seen as 
an algebraic equation in
$M^{(2)}.$ We introduce the following notation,
\begin{equation}
\Delta(S,t,u,u_S,u_t,u_{SS}, u_{St}, u_{tt})=u_t+\frac{\sigma^2
  S^2}2\frac{u_{SS}}{(1-\rho S u_{SS})^2}. \label{delta}
\end{equation}
Equation (\ref{urav}) is then equivalent to the relation
\begin{equation}
\Delta(w)=0, ~~ w \in M^{(2)}.
\end{equation}
We identify this algebraic equation with its solution
manifold $L_{\Delta}$ defined by
\begin{equation}
L_{\Delta}=\{ w \in M^{(2)} | \Delta(w)=0 \} \subset M^{(2)}.
\end{equation}

We consider an action of Lie-point groups on our differential equation
and its solutions.  We are interested in the group
$\rm{Diff}( M^{(2)})$ compatible with the contact structure of $
M^{(2)}.$ We denote the corresponding algebra  
by
$\mathcal Diff ( M^{(2)}).$
The symmetry group $G_\Delta$ of $\Delta$ is defined by
\begin{equation}
G_\Delta=\{ g \in \rm{Diff}( M^{(2)})|~~ g: ~~L_{\Delta} \to L_{\Delta}\}.
\end{equation}

\begin{theorem}
The differential equation (\ref{urav})
  admits a nontrivial four dimensional Lie algebra spanned by generators
$$ V_1 =  \frac{\partial}{\partial t}, ~~
V_2 =  S\frac{\partial}{\partial u},~~
V_3 = \frac{\partial}{\partial u}, ~~
V_4 = S \frac{\partial}{\partial S}+ u \frac{\partial}{\partial u}.$$ 
\end{theorem}
{\bf Proof}.
Let us consider a Lie-point vector field on $M,$ whose elements
are represented by
\begin{equation}
V= \xi (S,t,u) \frac{\partial}{\partial S} + \tau (S,t,u)
\frac{\partial}{\partial t} + \phi(S,t,u) \frac{\partial}{\partial u},
\end{equation}
where $\xi (S,t,u),\tau (S,t,u)$ and $\phi(S,t,u)$ are smooth
functions of their arguments, $V \in {\mathcal Diff} (M) $. Assume  there
exists an infinitesimal generator of an action $g \in G_{\Delta}. $ The
infinitesimal generators of these transformations form an algebra
${\mathcal Diff}_{\Delta} ( M).$ A Lie group of transformations acting
on the base space $M$ induces transformations on $M^{(2)}.$ The
corresponding algebra ${\mathcal Diff}_{\Delta} ( M^{(2)})$ will be
composed of the vector fields
\begin{eqnarray}
pr^{(2)} V &=& \xi (S,t,u) \frac{\partial}{\partial S} + \tau (S,t,u)
\frac{\partial}{\partial t} + \phi(S,t,u) \frac{\partial}{\partial u} \nonumber \\
&+&\phi^S(S,t,u) \frac{\partial}{\partial u_S}+\phi^t(S,t,u) \frac{\partial}{\partial u_t}  \label{prol2} \\
&+&\phi^{SS}(S,t,u) \frac{\partial}{\partial u_{SS}}+\phi^{St}(S,t,u) \frac{\partial}{\partial u_{St}}
+\phi^{tt}(S,t,u) \frac{\partial}{\partial u_{tt}} ,\nonumber
\end{eqnarray}
where $pr^{(2)} V$ is the second prolongation of the vector filed $V$.
Here the smooth functions $\phi^S(S,t,u)$, $\phi^t(S,t,u)$,
$\phi^{SS}(S,t,u)$, $\phi^{St}(S,t,u)$ and $\phi^{tt}(S,t,u)$ are
uniquely defined by the functions $ \xi (S,t,u), \tau (S,t,u)$ and
$\phi(S,t,u)$ using the prolongation procedure (see
\cite{Ovsiannikov}, \cite{Olver}, \cite{Stephani}, \cite{Gaeta},
\cite{Ibragimov}).

For our calculations we will use the explicit form of the coefficients
$\phi^t(S,t,u)$ and $\phi^{SS}(S,t,u)$ only because of the special
structure of equation (\ref{urav}).  The coefficient
$\phi^t(S,t,u)$ can be defined by the formula
\begin{equation}
\phi^t(S,t,u)= \phi_t+ u_t \phi_u -u_S \xi_t - u_S u_t \xi_u - u_t \tau_t - (u_t)^2 \tau_u \label{fit}
\end{equation}
and the coefficient $\phi^{SS}(S,t,u)$ by the expression
\begin{eqnarray}
\phi^{SS}(S,t,u)&=&\phi_{SS} + 2 u_S \phi_{S u} + u_{SS} \phi_u \label{fiss} \\
&+&(u_S)^2 \phi_{uu} - 2 u_{SS} \xi_S - u_S \xi_{SS}- 2 (u_S)^2 \xi_{Su}\nonumber \\
&-&3 u_S u_{SS} \xi_u -(u_S)^3 \xi_{uu} -2 u_{St} \tau_S -u_t \tau_{SS}\nonumber \\
&-&2 u_S u_t \tau_{Su}-(u_t u_{SS}+2 u_S u_{St}) \tau_{u} - (u_S)^2 u_t \tau_{uu} \nonumber ,
\end{eqnarray}
where the subscripts of $\xi, \tau , \phi$ denotes corresponding
partial derivatives.  The symmetry algebra ${\mathcal Diff}_{\Delta} (
M^{(2)})$ of the second order differential equation $\Delta=0$ can be
found as a solution of the determining equation
\begin{equation}
pr^{(2)} V(\Delta)=0 ~(mod(\Delta =0)), \label{algebradef}
\end{equation}
i.e., the equation (\ref{algebradef}) should be satisfied on the
solution manifold $L_{\Delta}$.  It is easy to prove that equation
(\ref{algebradef}) has the following solutions,
\begin{eqnarray}
V_1 &=& S \frac{\partial}{\partial S} +u \frac{\partial}{\partial u}, ~~~~
V_2 =  \frac{\partial}{\partial t}, \label{al4}\\
V_3 &=&  S\frac{\partial}{\partial u},~~~~~~~~~~~~~
V_4 = \frac{\partial}{\partial u}, \nonumber
\end{eqnarray}
where $V_i \in {\mathcal Diff}_{\Delta}(M), i=1,2,3,4.$
The commutative relations are
\begin{eqnarray}
&&~~[V_1,V_2]=[V_1,V_3]=[V_2,V_3]=[V_2,V_4]=[V_3,V_4]=0,\nonumber\\&&~~[V_1,V_4]=-V_4.\label{comut}
\end{eqnarray}
The vector fields $V_i, ~i=1,2,3,4$ span a four dimensional solvable
Lie algebra.  
\qed

An element of the algebra ${\mathcal Diff}_{\Delta} ( M)$ can be represented
 as a linear combination of the vector fields given by formulas (\ref{al4}) 
\begin{equation} \nonumber
V=a_1 V_1 +a_2 V_
2 +a_3 V_3 +a_4 V_4= \xi_a (S,t,u) \frac{\partial}{\partial S} + \tau_a (S,t,u)
\frac{\partial}{\partial t} + \phi_a (S,t,u) \frac{\partial}{\partial u},
\end{equation}
where
\begin{equation}\nonumber
\xi_a (S,t,u)= a_1 S,~\tau_a (S,t,u)=a_2,~\phi_a (S,t,u) =a_1 u + a_3 S + a_4
\end{equation}
with arbitrary constants $a_1,a_2,a_3,a_4$. 

Every element $V$ of the algebra ${\mathcal Diff}_{\Delta} ( M)$ is 
an infinitesimal generator of an action $g \in G_{\Delta}$.
Using the Lie equations we prove the following theorem.
\begin{theorem} 
\label{symteor}
The action of the symmetry group $G_{\Delta}$ of (\ref{urav})
 is given by (\ref{str})--(\ref{utro}). 
\end{theorem}

{\bf Proof}.
To find the transformations of the Lie group $G_\Delta$ associated with the 
generators (\ref{al4}) we just integrate the system
of ordinary differential equations, the so-called Lie equations,
\begin{eqnarray}\label{symsys1}
\frac{d{\tilde S}}{d\epsilon}=\xi_a ({\tilde S},{\tilde
t},u),~\frac{d{\tilde t}}{d\epsilon}=\tau_a ({\tilde S},{\tilde
t},{\tilde u}), ~\frac{d{\tilde u}}{d\epsilon}=\phi_a ({\tilde
S},{\tilde t},{\tilde u}),
\end{eqnarray}
with initial conditions
\begin{equation}\label{gransys1}
{\tilde S}|_{\epsilon=0}=S,~{\tilde t}|_{\epsilon=0}=t,~{\tilde u}|_{\epsilon=0}=u,
\end{equation}
where $\epsilon$ is the group parameter.
Here the variables ${\tilde S},{\tilde t}$ and ${\tilde u}$
denote the values $S,t,u$ after a symmetry transformation. The
solutions to the system of ordinary differential equations
(\ref{symsys1}) with initial conditions (\ref{gransys1}) have the
form
\begin{eqnarray}
{\tilde S}&=&S e^{a_1 \epsilon},~~ \epsilon \in (- \infty,\infty ),\label{str}\\
{\tilde t}&=&t+ a_2 \epsilon,\label{ttr}\\
{\tilde u}&=&u e^{a_1\epsilon} + a_3 S \epsilon e^{a_1 \epsilon} +
\frac{a_4}{a_1}(e^{a_1 \epsilon} -1), ~a_1 \ne 0 \label{utrne}\\
{\tilde u}&=&u +a_3 S \epsilon +a_4 \epsilon,~ a_1=0.\label{utro}
\end{eqnarray}
The equations (\ref{str})--(\ref{utro}) 
represent the action of the four parametric symmetry group $G_\Delta$.
\qed

 We will use this
symmetry group to construct invariant solutions to equation
(\ref{urav}). In detail the method of construction of invariant solutions  
is given in the book \cite{Ovsiannikov} and in the third chapter 
of the book \cite{Olver}. A lot of examples are given in the books
 \cite{Stephani}, \cite{Gaeta}, \cite{Ibragimov}.

To obtain the invariants of the symmetry group $G_\Delta$
 we can use a shortcut because of the 
very simple structure of the  Lie algebra found.

We exclude $\epsilon$ from the equations
 (\ref{str})--(\ref{utro}).
 Two
functionally  independent invariants can be taken in the form
\begin{eqnarray}
inv_1&=&a_1 t -a_2 \ln S,\label{invar1}\\
inv_2&=&a_1 \frac
{u}{S}  -a_3 \ln{S}+\frac{a_4}{S}, S>0.\label{invar2}
\end{eqnarray}

The functions (\ref{invar1})--(\ref{invar2}) are not defined at the point $S=0$
and, although the model equation (\ref{urav}) is defined at that 
point, we will exclude  $S=0$ in all further investigations.

We remark that the form of these invariants is not unique. Each function
of invariants (\ref{invar1}), (\ref{invar2}) will be an invariant.
Especially we can multiply each of the invariants by a constant
because any constant is a trivial invariant of the group $G_\Delta$.
But it is possible to obtain just two nontrivial functionally
independent invariants which we take in the form (\ref{invar1}),
(\ref{invar2}). The invariants can be used as new independent and 
dependent variables.

\section{Scaling variables}\label{scava}

Using the symmetry group $G_\Delta$ found in the preceding section
we reduce equation (\ref{urav}) to an ordinary differential
equation and define families of invariant solutions.
\begin{theorem} 
Up to the group transformations given by (\ref{str})-(\ref{utro}) 
all nontrivial Lie invariant solutions to equation (\ref{urav}) 
depend on the scaling variables $z, v(z)$ and the relations
\begin{eqnarray}
z&=&  \ln S-\delta t,~~ \delta \ne 0,\label{trans_z}\\
u(S,t)&=& - S v(z), \label{trans_v}
\end{eqnarray}
where $\delta$ is an arbitrary constant, hold.
\end{theorem}
{\bf Proof}.
We can reduce the partial differential equation (\ref{urav}) to
an ordinary differential equation for the function $v(z)$
if we change the variables $u,S,t$ for
$z=\phi(S,t,u)$ and $v=\psi(S,t,u)$. This substitution leads to 
invariant solutions to equation (\ref{urav}) if $\phi(S,t,u)$ and
$\psi(S,t,u)$ are some invariants of the symmetry group $G_{\Delta}$.
In the previous section we found just two invariants, hence all 
invariant solutions except for trivial ones will arise after the substitutions
(\ref{trans_z})--(\ref{trans_v}). The trivial solutions we can obtain if
we assume  $u={\rm const.}$, $u=u(t)$ and $u=u(S)$.

We remark that we take as a new independent variable the first
invariant (\ref{invar1}) of the symmetry group $G_\Delta$ 
and as the dependent variable the nontrivial part of the second invariant,
this allows us to simplify the calculations. In this way we do not 
lose any solutions because 
the found invariant solutions can be 
later transformed by the 
rule of thumb given by (\ref{utrne})--(\ref{utro}).
\qed

 The equation for
the function $v(z)$ has the form
\begin{eqnarray}
\label{uravz}
v_z\,(1+\rho  \,(v_z+ v_{zz}))^2 - \frac{\sigma^2 }{2\delta} \,(v_z+ v_{zz})=
0.
\end{eqnarray}

The Lie group of symmetries  for this equation can be found in the same way as  described
in previous Section \label{liesec} for equation (\ref{urav}).
\begin{theorem} \label{ordsym} \cite{Lie}
The equation (\ref{uravz}) admits a two dimensional Abelian Lie algebra spanned by two generators
\begin{equation}
\label{gen}
U_1= \frac{\partial}{\partial z},~~ U_2= \frac{\partial}{\partial v}.
\end{equation}
\end{theorem}
{\bf Proof}.
This theorem was proved in a more general case  by S.Lie in (\cite{Lie}). Also
it can be verified by a straightforward calculation.
\qed

Equation (\ref{uravz}) allows a two-dimensional Lie group associated with the
Lie algebra spanned by the generators (\ref{gen}).  As a consequence equation
(\ref{uravz}) is completely integrable. Hence the most general form of the
solution of (\ref{uravz}) is a two parametric family of congruent curves. To
obtain a two parametric family of solutions to equation (\ref{uravz}) we can
subsequently use the two generators (\ref{gen}) in arbitrary order.  Both
ways will lead to the same family of solutions independent on the order.  To obtain a solution we must perform two integrations and this
procedure is not always possible in closed form. However, in view  of the theorem
\ref{ordsym} we do not have any other possibility to solve equation
(\ref{uravz}) in a more convenient
way.

In the next Section we put constraints on the
constant $\delta$ in (\ref{trans_z}) in order to integrate the arising 
equations in an exact form.
Consequently we restricted ourselves and do not obtain the most
general form for the family of solutions.

\section{Families of invariant solutions}\label{fainv}

\begin{theorem}\label{thurav}
The equation (\ref{uravz})
can be reduced by the substitution $v(z)_z=y(z)$ to the set of equations 
\begin{eqnarray}
 y(z)&=&0,~~~~y(z)=\frac{1}{\rho}\left(-1 \pm \sqrt{\frac{\sigma^2}
{2 \delta}}\right), \label{trivisol}\\
\frac{dy}{dz}&=&-\frac{1}{ y}
\left(\left(y^2+\frac{y}{\rho} -\frac{\sigma^2}{4\rho^2 \delta}
\right) \pm \frac{1}{\rho} \sqrt{\frac{\sigma^2}{2 \rho \delta}} 
\sqrt{\frac{\sigma^2}{8 \rho \delta}-y}  \right), ~ y\ne 0,\label{secst}
\end{eqnarray}
where $\delta$ is an arbitrary constant.\\
The complete set of solutions to equation (\ref{uravz}) coincides with the union of solutions to these equations.
\end{theorem} 
{\bf Proof}.
First we look for the solutions of the type $v(z)_z={\rm const.}$
From  straightforward calculations we obtain that equations (\ref{trivisol}) hold. The corresponding solutions to equation (\ref{uravz}) have the form
\begin{equation}\label{consolv}
v(z)=c_1,
\end{equation}
where $c_1$ is an arbitrary constant,
and
\begin{equation}  \label{firv}
v(z)=-\frac{1}{ \rho} \left(1 \pm \sqrt{\frac{\sigma^2}{2 \delta}} \right)z + {\rm const}.
\end{equation}

We assume now that $y(z) \ne 
{\rm const.}$, i.e., $v_z(z) \ne {\rm const.} $ and
 use  the operator $U_2$, (\ref{gen}), first to introduce the new dependent variable
$y(z)=v_z(z)$ in equation (\ref{uravz}). 

We obtain a first order differential
equation for the function $y(z)$,
\begin{eqnarray}
\label{urav_yz}
y_z^2+2 \frac{y_z}{y} \left(y^2+\frac{y}{\rho} -\frac{\sigma^2}{4\rho^2 \delta}\right)+ \left(y^2+\frac{2}{\rho}\,y+\frac{2 \delta-\sigma^2}{2 \rho^2 \delta} \right) =0, ~~y \ne 0.
\end{eqnarray}

The equation (\ref{urav_yz}) is quadratic in the highest derivative and it
can be represented as a product of two differential equations (\ref{secst}).

We reduced equation (\ref{urav_yz}) to a product
of two equations (\ref{secst}) and  in
this way we could have lost some of the solutions.
Let us now study the
discriminant curve for equation (\ref{urav_yz}). 
We denote by $F(y_z,y,z)$ the
left hand side of equation (\ref{urav_yz}), i.e.,
\begin{equation}
F(y_z,y,z)=y_z^2+2 \frac{y_z}{y} \left(y^2+\frac{y}{\rho} -\frac{\sigma^2}{4\rho^2 \delta}\right)+ \left(y^2+\frac{2}{\rho}\,y+\frac{2 \delta-\sigma^2}{2 \rho^2 \delta} \right).
\end{equation}
The discriminant curve is a set of points fulfilling the conditions,
\begin{eqnarray}
F(y_z,y,z)=0, \label{ur}\\
\frac{\partial F(y_z,y,z)}{\partial y_z} =0. \label{os}
\end{eqnarray}
Along this curve the conditions of the theorem on an implicit
function are not satisfied and  in these points
 the obtained solutions may be not unique.
It is easy to prove that the system of equations (\ref{os})--(\ref{ur}) 
has a unique solution,
\begin{equation}\label{exep}
y_{excep}(z) = \frac{1}{\rho},
\end{equation}
for the special value of the constant $\delta$ 
 \begin{equation}\delta= \sigma^2/8\label{delta}
\end{equation}only.
The corresponding solution of equation (\ref{uravz})
has the form
\begin{equation} \label{spec}
v_{excep}(z) = \frac{ z}{ \rho} +{\rm const}
\end{equation}
and it coincides with one of the solutions (\ref{firv}) for 
 $\delta= \sigma^2/8.$
\qed

\begin{theorem}
The explicit invariant solutions to equation (\ref{urav}), defined on the
region $S>0, t \in [0,T], ~T>0$ are given by (\ref{consolv}),(\ref{firu}) and
(\ref{scalsol}). Other solutions of this type can be obtained using the
transformations of the symmetry group $G_\Delta$ represented by
(\ref{utrne})--(\ref{utro}).
\end{theorem}
{\bf Proof}.
To obtain the invariant solutions we should solve the equations listed
in the theorem \ref{thurav}. It is trivial to solve the first two of them.

The
relations (\ref{trivisol}) have
the following solutions
\begin{equation}\label{consolv}
u(S,t)= S c_1,
\end{equation}
and
\begin{equation} \label{firu}
u(S,t)=\rho^{-1} \left(1 \pm \sqrt{\frac{\sigma^2}{2 \delta}}\right) (S \ln{S} - \delta S t) + S d_0,
\end{equation}
where $\delta$ and $d_0$ are arbitrary constants.

To integrate the equations (\ref{urav_yz})  we use the second operator $U_1$, (\ref{gen}), and separate variables
\begin{eqnarray}
 \int{y \left(\left(y^2+\frac{y}{\rho} -\frac{\sigma^2}{4\rho^2 \delta}
\right) \pm \frac{1}{\rho} \sqrt{\frac{\sigma^2}{2 \rho \delta}} 
\sqrt{\frac{\sigma^2}{8 \rho \delta}-y}  \right)^{-1}}
=-z+{\rm const.} \nonumber 
\end{eqnarray}
We denote the integral on the left hand side by 
\begin{eqnarray}
\label{int_y}
I(y)=\int{y \left(\left(y^2+\frac{y}{\rho} -\frac{\sigma^2}{4\rho^2 \delta}
\right) \pm \frac{1}{\rho} \sqrt{\frac{\sigma^2}{2 \rho \delta}} 
\sqrt{\frac{\sigma^2}{8 \rho \delta}-y}  \right)^{-1}}.
\nonumber \end{eqnarray}
Straightforward calculations lead  to 
the following form for the function $I(y)$,
\begin{eqnarray}
\label{loesung_xi}
I(y)=\frac{2 \delta}{2 \delta - \sigma^2}\left(-\frac{\sigma^2}{2 \delta}
\,\ln\left(\xi\mp\frac{1}{2}\sqrt{\frac{\sigma^2}{2 \rho \delta}}
\right)\right.&+&\left. 
\frac{4 \delta - \sigma^2}{4 \delta}
\ln\left(\xi^2\pm\sqrt{\frac{\sigma^2}{2 \rho \delta}}
\,\xi-\frac{8 \delta - \sigma^2}{8 \rho \delta}
\right)\right.\nonumber\\
&\mp&\left.\sqrt{\frac{\sigma^2}{2 \delta}}
{\rm arctahn}\left( \sqrt{\rho}\left(\xi\pm\frac{1}{2}\sqrt{
\frac{\sigma^2}{2 \rho \delta}}\right)\right)\right),
\nonumber \end{eqnarray}
where the variables $\xi$ and $y$ are connected by
\begin{equation}
\label{transform_xi}
\frac{\sigma^2}{8 \rho \sigma}-y =\xi^2.
\end{equation}

Now let us chose 
$ \delta=\sigma^2/8 $, i.e.,  in a way that substitution (\ref{trans_z})
takes the form
\begin{equation}
z=\ln S - \frac{\sigma^2}{8}t .
\end{equation}
 We obtain an explicit representation for the
function $y(z)$ which solves equation (\ref{urav_yz}).  The
solutions are given by
\begin{equation} \label{resh}
y(z)=-\frac{1}{\rho} \left(1+\frac{2^{4/3} e^{z} } {\left(\left(m+ e_1 \sqrt{m^2+4e^{3z/2}}\right)^{4}\right)^{1/3}}+ \frac
{\left(\left(m+ e_1 \sqrt{m^2+4e^{3z/2}}\right)^{4}\right)^{1/3}}
{2^{4/3} e^{z}}\right)
\end{equation}
with an arbitrary constant $m$ and $e_1 = \pm 1$.

Thereafter we integrate the equation $v_z=y(z)$ and obtain a family of
solutions to equation (\ref{uravz}),
\begin{eqnarray}
v(z)=-\frac{1}{\rho}\left(z-2^{-4/3} e^{-z} \left(\left( m+ \epsilon_1\sqrt{m^2 + 4 e^{\frac{3}{2 } z}}\right)^{4}\right)^{1/3} \right. \nonumber\\
\left. -2^{-4/3} e^{-z}
\left( \left(- m+ \epsilon_1\sqrt{m^2 + 4 e^{\frac{3}{2 } z}}\right)^{4}\right)^{1/3} \right. \nonumber\\
\left. - \ln{\left(\left( m+\epsilon_2\sqrt{m^2 + 4 e^{\frac{3}{2} z}}\right)^{1/3}-
            \left( -m+\epsilon_2\sqrt{m^2 + 4 e^{ \frac{3}{2}z}}\right)^{1/3} \right)^4}+d \right), \nonumber
\end{eqnarray}
where $d$ and $m\ne 0$ are arbitrary constants. The case $m=0$
corresponds to the solution (\ref{firv}).
The parameters $\epsilon_1$ and $\epsilon_2$
take values $\epsilon_1 = \pm 1$ , $\epsilon_2 =\pm 1$ and can be
chosen independently. The solutions do not depend on the value of
$\epsilon_1$.

The corresponding family of solutions to equation (\ref{urav}) will
take the form
\begin{eqnarray}
 u(S,t)&=&\rho^{-1} S \ln{S} - \frac{\sigma^2}{8 \rho}   S t  \nonumber \\
&
&- 2^{-4/3} \rho^{-1} \exp { \left( \frac{\sigma^2 t}{8} \right)}
\left(\left(m+\epsilon_1 \sqrt{m^2 +4 S^{3/2} \exp \left(-\frac{3\sigma^2 t}{16}\right)}\right)^{4} \right)^{1/3}
 \nonumber\\
&&- 2^{-4/3} \rho^{-1} \exp \left(\frac{\sigma^2 t}{8}
                         \right)
\left(
\left(- m+\epsilon_1  \sqrt{m^2 +4 S^{3/2} \exp \left(-\frac{3\sigma^2 t}{16}
                                              \right)}
\right)^4
\right)^{1/3} \nonumber \\
&&- \rho^{-1} S \ln \left( \left(m+\epsilon_2\sqrt{m^2 + 4 S^{3/2} \exp \left(-\frac{3 \sigma^2 t}{16}\right)}\right)^{1/3}
\right.
\label{scalsol}  \\
&& \left.~~~~~~~~~~~~ - \left(- m+\epsilon_2\sqrt{ m^2 + 4 S^{3/2} \exp\left(-\frac{3\sigma^2 t}{16} \right) }\right)^{1/3} \right)^4  + S d_1 +d_2, \nonumber
\end{eqnarray}
where $d_1,d_2$ are arbitrary constants. In formula (\ref{scalsol}) we
assume that the arbitrary parameter $m$ is non equal to zero. In case
$m=0$ this solution can be reduced to one of the solutions
(\ref{firu}).

Let us compare this solution with the solutions to equation (\ref{urav_yz})
which were obtained in the case $\delta=\sigma^2/8$.
The functions $y(z)$ in the family (\ref{resh}) are even functions of the constant $m$.
For $m=0$ we obtain
\begin{equation}\label{exep2}
y(z) = -\frac{3}{\rho}.
\end{equation}
The solution (\ref{exep2}) leads to the described solutions (\ref{firv}) and
(\ref{firu}) with upper sign and with  $\delta = \sigma^2/2$.

\qed

\section{Properties of invariant solutions}\label{prop5}

The solutions (\ref{scalsol}) depend on the parameter $\rho$ in a very
simple way: all solutions of this
family have the factor $1/ \rho$ in front of the whole expression.
This parameter, which is the measure of the influence of the
large trader on the market, is a constant $0 \le \rho $ and cannot
be equal to zero for the large trader.   This
means that each solution of this family does completely blow up at $\rho
\to 0$. Consequently these solutions have no linear analogies.
 
If we denote by $\tilde u(S,t)=\rho u(S,t) \rho$ we obtain for the function
$\tilde u(S,t)$ following equation
\begin{eqnarray} \label{uravohner}
\tilde u_t+\frac{\sigma^2 S^2}2\frac{\tilde u_{SS}}{(1- S \tilde u_{SS})^2}=0.
\end{eqnarray}

This means that the solutions (\ref{scalsol}) multiplied by $\rho$ are solutions to
 equation (\ref{uravohner}). If we obtain any solutions to (\ref{uravohner})
 for any fixed boundary conditions, we obtain the corresponding solutions to
 (\ref{scalsol}) with boundary conditions divided by $\rho$ if we divide the
 found solutions by $\rho$ as well. In other words, the solutions (\ref{scalsol})
 strongly reflect to the nonlinearity in equation (\ref{urav}).

\begin{figure}[ht]
\vspace{0.5cm}
\begin{minipage}[t]{14cm}
\includegraphics[width=14cm]{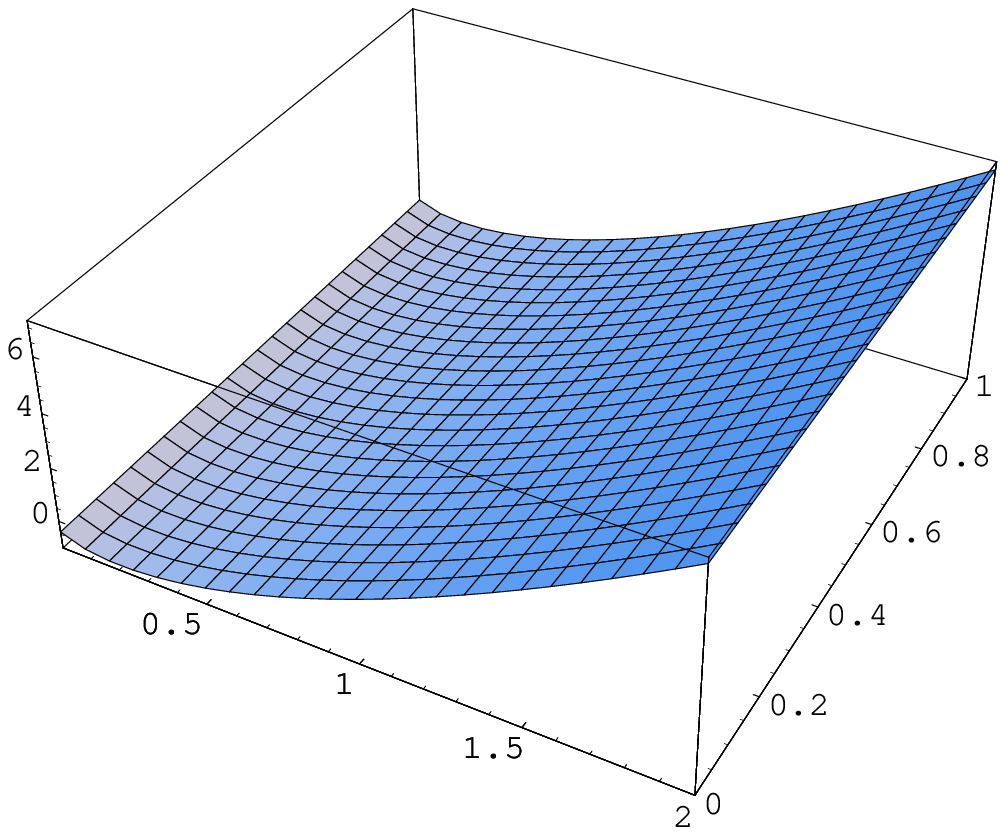}
\caption{\label{solfi1} Plot of the solution $u(S,t) $ (\ref{scalsol}) with $S \in (0,2]$, $t \in [0,1]$ and parameters $\sigma=0.35,~m=0.5, \rho=0.1$, $d_1=d_2=0.$}

\end{minipage} 
\end{figure}

Let us study the analytical properties of the solution  (\ref{scalsol})
and the corresponding payoff.
In Figure \ref{solfi1} we represent graphically  the solution $u(S,t)$
(\ref{scalsol}) for small values of the variables $S,t$.

\begin{figure}[ht]
\begin{minipage}[t]{14.cm}
\includegraphics[width=14.cm]{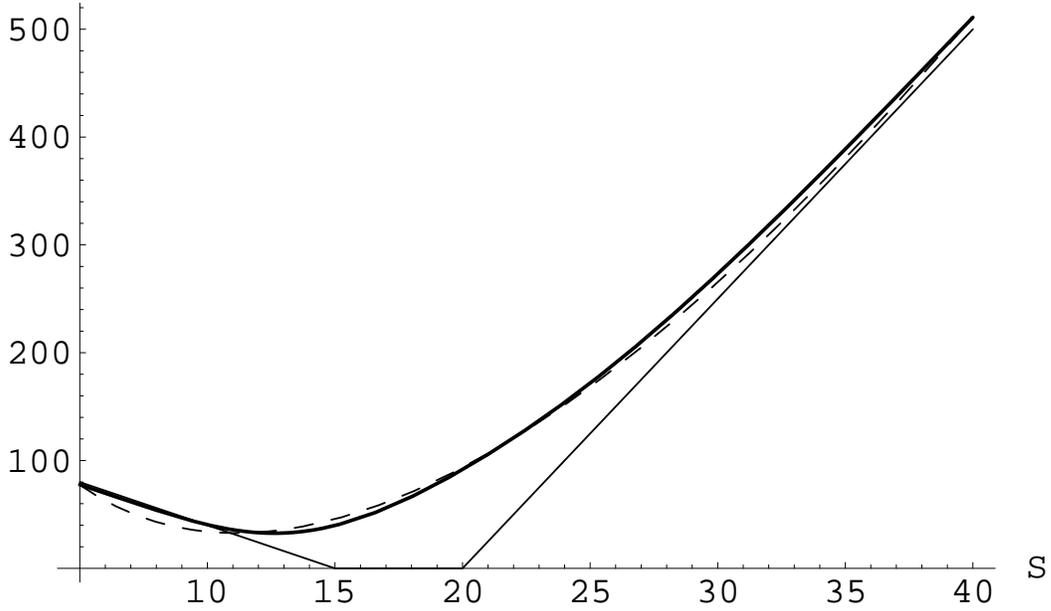}
\caption{\label{solstr1} 
Plot of the explicit solution $u(S,t)$, (\ref{scalsol}),(dashed line) with 
 parameters $\rho = 0.05,~ m=1338.0,~d_1=140.0,~ d_2=295139,~t=0$
compared with the solution
$u(S,0)$ (solid line) of the linear Black--Scholes model (\ref{BS}) with the payoff $u_{strangle}(S,T)$ (thin solid line). The parameters for the strangle are 
$r=0.02$, $~ \sigma=0.25$, $E_P=15.0$, $~E_C=20.0~$ and $T=1.0$.}
\end{minipage}
\end{figure}

Let us represent the payoff of a strangle as a sum of $K_P$ European puts with
an exercise price $E_P$ and $K_C$ European calls with an exercise price $E_C$
which have the same expiry date $T$.
We can choose in an appropriate way the parameters $m,d_1,d_2$
of the explicit solution (\ref{scalsol}) such
that this solution approximates the payoff of 
a strangle  $u_{strangle}(S,T)$, 
\begin{equation} \label{stranglepayoff}
u_{strangle}(S,T)=K_P \max{\left(E_P-S,0\right)} + K_C \max{\left(S-E_C,0\right)}, ~~ E_P<E_C.
\end{equation}
This is shown in Figure \ref{solstr1}.

Let us now investigate the asymptotic properties of solutions (\ref{scalsol})
for $S \to 0$ and for $S \to \infty$.
The
asymptotic behaviour of the function (\ref{scalsol}) for $S \to 0$ can be
described as follows,
\begin{eqnarray}
\rho u(S,t) &\sim& (m^{4})^{1/3} \exp\left(\frac{\sigma^2}{8} t\right) + S \ln S \label{asynul}\\
&-& S \left( \frac{\sigma ^2 t}{8} + \ln ((2 m)^{4})^{1/3} \right)
- S^{3/2}\frac{4 \exp{\left(-\frac{\sigma ^2 t}{16}\right)}}{3 (m^{2})^{1/3}}  +O(S^{5/2}) ,~~ S \to 0.\nonumber
\end{eqnarray}
The main term in formula (\ref{asynul}) depends on the time and on the
constant $m.$ We can choose $m$ to model payoff properties.  From this
decomposition it follows immediately that for all solutions from this
family the denominator in equation (\ref{urav}) vanishes in the point
$S=0$.  In order to avoid this singularity we exclude the point $S=0$ from
the intervals where the numerical investigations are done.  \\[5pt]

The main term of the asymptotic expansion of $u(S,t)$ for $ S \to \infty$,
\begin{eqnarray}
\rho ~u(S,t) &\sim& 3 S \ln S - S \left( \frac{3 \sigma ^2 t}{8}+
4 \ln \left(\frac{2^{1/3} m }{3} \right) +2
\right) \nonumber\\
&-& S^{-1/2}\left(\frac{2}{3}\right)^3 \exp{\left(\frac{3\sigma ^2 t}{16}\right)} m^2
 +O(S^{-5/4}) ,~~ S \to \infty, \label{asybes}
\end{eqnarray}
is equal to $3 S \ln S.$ This term is independent of any integration
constant or time.  Hence all solutions in this family have the same
asymptotic behaviour for $S \to \infty$.

\begin{figure}[ht]
\centerline{The sensitivity parameters $\Theta$ and {\it Vega} for solutions (\ref{scalsol}).}
\vspace{0.5cm}
\begin{minipage}[t]{7cm}
\includegraphics[width=7cm]{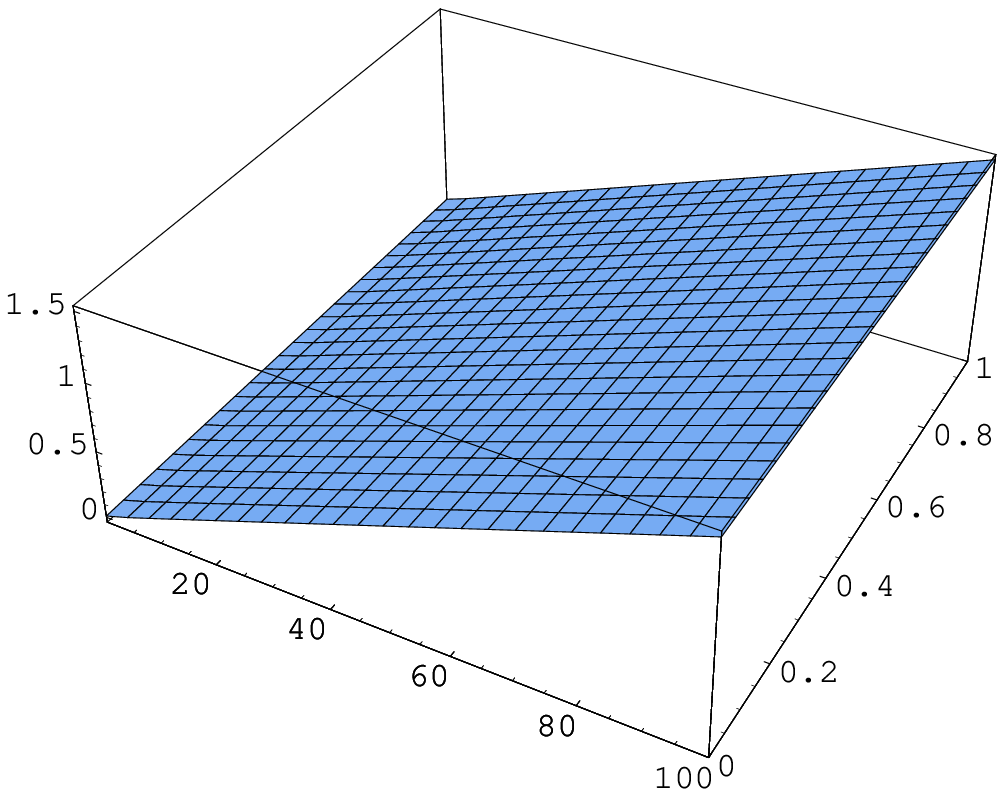}
\caption{\label{theta} Plot of
$\rho \Theta =-\rho \frac{\partial u(S,t)}{\partial t}$ with $S \in (0,100]$, $t \in [0,1]$ 
and parameters $\sigma=0.2,~m=-1.7$, $d_1=d_2=0.$ }
\end{minipage} \hfill
\begin{minipage}[t]{7.cm}
\includegraphics[width=7.cm]{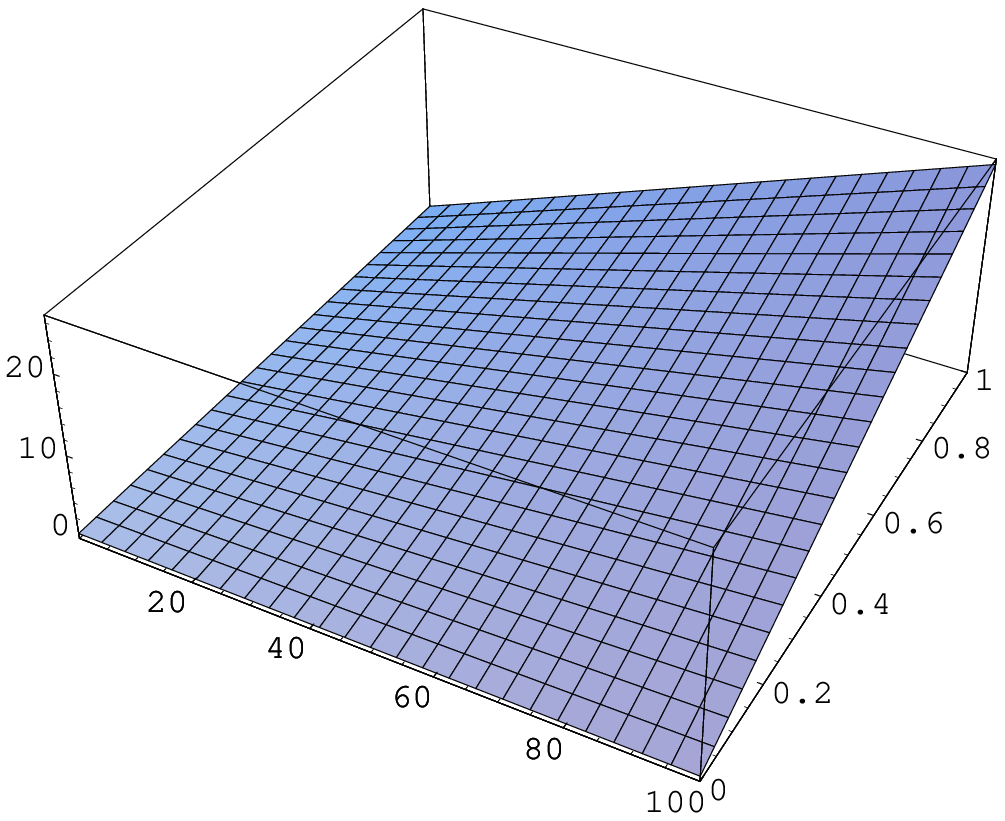}
\caption{\label{vega} Plot of 
${\rho \it Vega} = -\rho \frac {\partial u(S,t)}{\partial \sigma} $ 
with $S \in (0,100]$, $t \in [0,1]$ and parameters $\sigma=0.35,~m=0.5$, $d_1=d_2=0.$}
\end{minipage}
\end{figure}

From the financial point of view it is important to study the
dependencies of the obtained solutions on different parameters,
for instance, on time, on volatility or on the price of the
underlying asset, etc. 
In this way we get information about the sensitivity of
our product  with respect to a change of one of these parameters.
Using the explicit formula for the solutions (\ref{scalsol}) it is
easy to represent these dependencies graphically, see Figures
\ref{delta} - \ref{vega}. The time dependence of the solutions 
(\ref{scalsol}) is very weak 
but still present as we can see on Figure \ref{vega}.

The obtained family of solutions (\ref{scalsol}) can be used 
as a benchmark for testing of numerical methods.
We suggest the following
procedure.
We use the solutions (\ref{scalsol}) with
boundary conditions which we can obtain just by fixing the time to test
numerical methods. These boundary conditions are smooth.  Then we take one of the
numerical methods and try to reproduce the analytical solution. In
this manner we can check on each time step the reached accuracy and
adjust the parameters of the grid and the numerical scheme.  Thereafter we
can be sure that for all smooth boundary conditions of the same type
as studied we obtain  numerical solutions with nearly the same
accuracy.  

Now if we apply this method to boundary conditions with worse properties we
can be sure that it works at least in the case of an approximation of these
boundary conditions by very close but smooth ones. 
  In case of
an European call option we have a continuous payoff function $u_1(S,T)$,
\begin{equation}
u_1(S,T)=\max( S-E, 0),  \label{one}
\end{equation}
where $E$ is the exercise price and $T$ is the expiry date, which is not
differentiable in $S=E.$ We can make it smooth by just replacing the
payoff $u_1(S,T)$ in the neighbourhood of the exercise price
$E$. Usually one takes as such smooth function a solution of the
linear Black--Scholes formula (\ref{BS}) for $t \sim T.$ Then we can
compare the results of numerical calculations in both cases. If they
do not have any significant difference we can use the same method also
in case of continuous boundary conditions and relax the condition of
smoothness.\\[4pt]

As a first example we take an explicit method for a numerical solution of
 equation (\ref{urav}).
 This method can be used to find numerically solutions to
 the linear Black--Scholes model
\begin{equation}\label{BS}
u_t +\frac{\sigma^2}{2}S^2 u_{SS} +r S u_S -r u=0,
\end{equation}
where $r$ is the interest rate.  It gives proper results for the
special relation between $\Delta S^2$ and $\Delta t$, where by $\Delta
S,\Delta t $ we denote correspondingly the mesh sizes of the
discretization of $S$ and $t$ intervals. We applied this method to the
nonlinear equation (\ref{urav}). We proved that in all studied cases
the explicit method diverges independently from the chosen relation
between $\Delta S^2$ and $\Delta t$.  It follows that the
explicit method is not reasonable in this nonlinear case.\\[4pt]

Another  way to solve equation (\ref{urav}) numerically is to
use the completely implicit method. For the linear Black--Scholes model (\ref{BS})
it gives proper results for arbitrary relations between $\Delta S^2$
and $\Delta t.$ For a nonlinear equation this method leads to a system
of nonlinear algebraic or transcendental equations.
An attempt to solve such a system can easily
exceed the possibilities of a modern computer due to the very fast with growing
 grid size. We used this method to reproduce the explicit
solutions (\ref{scalsol}) with appropriate accuracy.
Thereby we used equidistant grids with
$16,28$ and $42$ space nodes and with $15$ and $30$ time levels.
We reached the  relative accuracy of order of $0.2 \%$ .

Then we used this completely implicit method to calculate the value of
derivatives governed by equation (\ref{urav}) with usual payoff
functions.\\ 
Let us describe shortly the  system of difference
equations which we used. It was obtained by replacing  
 the derivatives in the  $t$ and $S$ directions in the following
way,
\begin{eqnarray}
\frac{\partial u}{\partial t}&=&\frac{u(S_i,t_{j+1})-u(S_i,t_{j})}{\tau}+O(\tau), \nonumber \\
\frac{\partial^2 u}{\partial
S^2}&=&\frac{u(S_{i+1},t_j)-2u(S_i,t_j)+u(S_{i-1},t_j)}{h^2}+O(h^2),
\end{eqnarray}
 
where $\tau= \Delta t$ is the time step and $h=\Delta S$ the 
 space step.
For each fixed $j$ we obtain a system of $N_S-1$ equations
\begin{eqnarray}\label{sysin}
&&\frac{u_{ij+1}}4\,\left(\frac{h^2}{S_i}-\rho (u_{i-1j}-2 u_{ij}+u_{i+1j})\right)^2-\frac{u_{ij}}4\,\left(\frac{h^2}{S_i}-\rho (u_{i-1j}-2 u_{ij}+u_{i+1j})\right)^2\nonumber\\
&&~~~~~~+\frac{\tau \sigma^2 h^2}{8}\,(u_{i-1j}-2 u_{ij}+u_{i+1j})
 =0,\qquad i=\overline{2,N_S},\,\,j=\overline{N_t,1},
\end{eqnarray}
for the internal points, where $N_t$ is the number of time layers and $N_S+1$ is the
number of grid nodes in space direction. In this case we used  the final
conditions, i.e., the knowledge of the values $u(S,T)$ and calculated the values for $u(S,t)$ backwards to $t=0$.
In the system (\ref{sysin}) the values on the layer $j+1$ are known and the
values on the layer $j$ are unknown functions.
On the boundaries $S_1$ and $S_{N_S+1}$  the values $u_{ij} $ are
defined for each fixed $j$ by the function $u_{bound}(S,t)$ in accordance with
 the used boundary conditions.
The complete system
  of difference equations has the form

\begin{eqnarray}
&&\frac{u_{2j+1}-u_{2j}}4\,\left(\frac{h^2}{S_2}-\rho (u_{bound}(S_1,t_j)-2 u_{2j}+u_{3j})\right)^2\nonumber\\
&&\phantom{\frac{u_{2j+1}-u_{2j}}4 \,\left(\frac{h^2}{S_2}\right.}+\frac{\tau \sigma^2 h^2}{8}\,(u_{bound}(S_1,t_j)-2 u_{2j}+u_{3j})
 =0,\label{i2}
\end{eqnarray}

\begin{eqnarray}
&&\frac{u_{ij+1}-u_{ij}}4\,\left(\frac{h^2}{S_i}-\rho (u_{i-1j}-2 u_{ij}+u_{i+1j})\right)^2\nonumber\\
&&\phantom{\frac{u_{ij+1}-u_{ij}}4}+\frac{\tau \sigma^2 h^2}{8}\,(u_{i-1j}-2 u_{ij}+u_{i+1j})
 =0,\quad i=\overline{3,N_s-1},\label{ins2}
\end{eqnarray}

\begin{eqnarray}
&&\frac{u_{N_sj+1}-u_{N_sj}}4\,\left(\frac{h^2}{S_{N_s}}-\rho (u_{N_s-1j}-2 u_{N_sj}+u_{bound}(S_{N_s+1},t_j))\right)^2\nonumber\\
&&\phantom {\frac{u_{N_sj+1}-u_{N_sj}}4 } +\frac{\tau \sigma^2 h^2}{8}\,(u_{N_s-1j}-2 u_{N_sj}+u_{bound}(S_{N_s+1j},t_j))
 =0 \label{ins}
\end{eqnarray}
with $j=\overline{N_t,1}$.

\begin{figure}
\begin{minipage}[t]{7cm}
\includegraphics[width=7.cm]{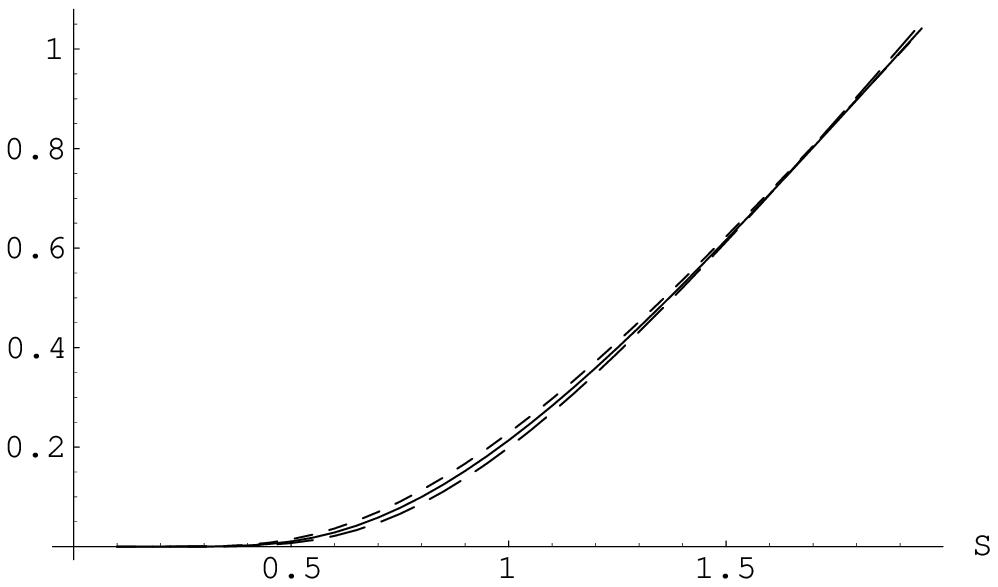}
\caption{\label{rhodep} Plot of the numerical solution  
for the hedge-cost $u_1(S,0)$, (\ref{urav}), with the payoff $u_{1}(S,T)$ (\ref{one})
 for various values of $\rho$. Compare the
solution $u_1(S,0)$  with $\rho=0.3$ (short dashed line), with $\rho=0.2$ (solid line) and with $\rho=0.1$ (long dashed line).}
\end{minipage}\hfill
\begin{minipage}[t]{7cm}
\includegraphics[width=7.cm]{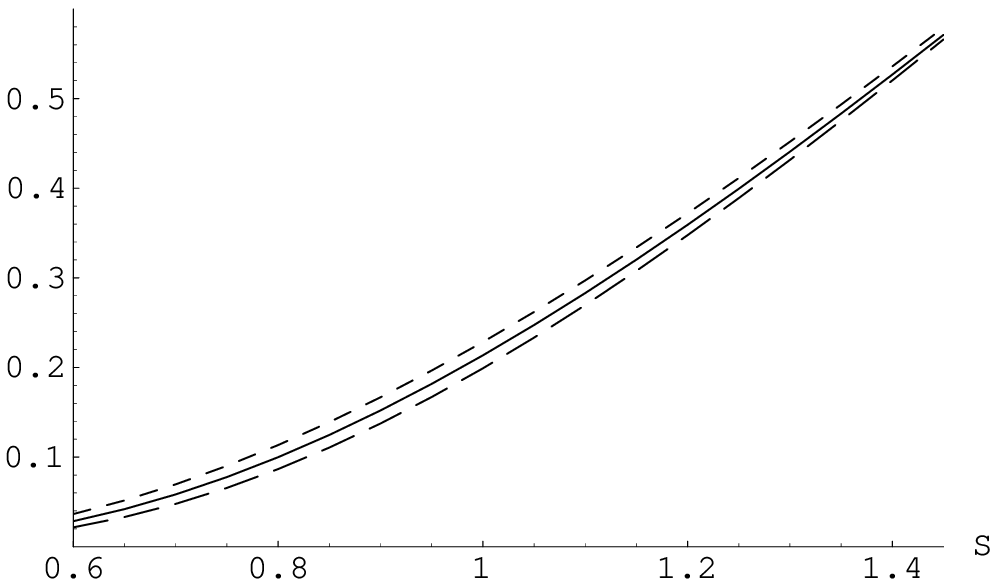}
\caption{\label{rhodep3} The part of the same curves as in 
Figure \ref{rhodep} in the neighborhood $S
  \sim E$.
The parameters of the European call are $\sigma=0.35$, $~T=0.9$, $~E=0.914$,  $ S=[0.1,2]$.
The parameters of the grid are $ ~h=0.05,~ \tau=0.05$, $~N_S=38$, $~ N_t=18.$
}
\end{minipage}
\end{figure}

\begin{figure}
\begin{minipage}[t]{14.cm}
\includegraphics[width=14.cm]{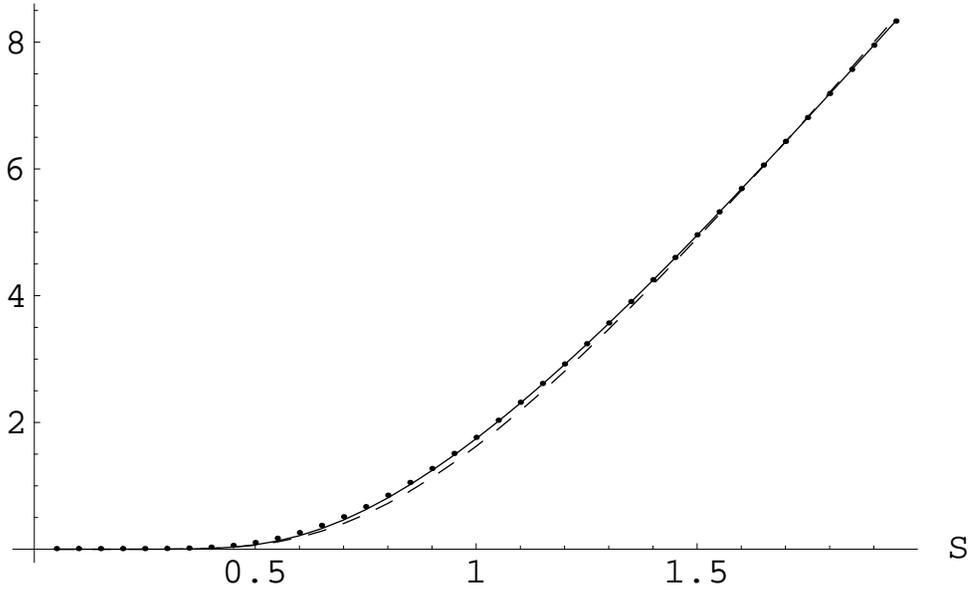}
\caption{\label{num} 
Plot of the numerical solution $u_8(S,0)$,  (\ref{urav}), (dots) 
 with the payoff $u_{8}(S,T)$, (\ref{ten}), associated with  8 European calls
compared with the solution
$u(S,0)$ (solid line) of the linear Black--Scholes model (\ref{BS}) with the same payoff, with the
sum of numerical solutions $ u_3(S,0)+ u_5(S,0)$ (dashed line) of
equation (\ref{urav}).
The parameters are $S \in [0.1,2.0]$, $t \in [0,T]$, $T=0.9$, $r=0.02$, 
$\sigma=0.35$, $~E=0.914$, $~\rho = 0.03$.
The parameters of the grid are $N_S=38, N_t=18$, $~ \tau=0.05$, $~ h=0.05.$}
\end{minipage}
\end{figure}

In the works \cite{Frey:perfect} and \cite{Frey:illicuidity} it was
proved that the hedge-cost of the claim $u(S,t) $ increases
monotonously with growing $\rho$, i.e., with growing influence of a
large trader. We prove this dependence numerically. We take as
boundary conditions $u_{bound}(S,T)=u_1(S,T)$ (\ref{one}), i.e., the
boundary conditions which correspond to one European call option. We
calculate the values $u(S,t=0)$ for various values of $\rho$. In
Figures \ref{rhodep} - \ref{rhodep3} we can see that with the growing
value $\rho$ the option hedge-cost also grows monotonically. It
completely corresponds to the functional behaviour obtained in the work
\cite{Frey:illicuidity}.\\[5pt] 
Now we  compare the option
hedge-cost predicted by the linear Black--Scholes model (\ref{BS}) and
by the nonlinear model (\ref{urav}).\\ 
At first we find numerically
the value of the hedge-cost for the derivative $u_3(S,0)$ defined by
equation (\ref{urav}) with the payoff $u_3(S,T)$ which corresponds to
 three European call options. The payoff function for $K$ European
call options is given by
\begin{equation}
u_{K}(S,T)=K \max(S- E, 0),  \label{ten}
\end{equation}
where we will use $K=3,5,8.$ Then we find numerically the value of the
 hedge-cost for the derivative $u_5(S,0)$ defined by equation
 (\ref{urav}) with the payoff $u_5(S,T)$ (\ref{ten}) which corresponds
 to  five European call options. The exercise price we take equal to
 $E=0.914$ in both cases.\\ 
Thereafter we calculated numerically the
 value of the hedge-cost of the derivative $u_{8}(S,t)$ with a payoff
 function which corresponds to the eight European call options with
 the same value $E=0.914$ as before and the same expiry date $T=0.9$
 and the same value $\rho=0.03$.  In these cases we use the grid with
 $N_S=38$, $N_t=18$, i.e., with $39$ nodes in the space direction and
 with $18$ time layers.  In the linear case it makes no sense to
 calculate once more the value for this derivative, we may just add
 the values $u_3(S,t)$ and $u_5(S,t)$ obtained in the former cases.
 However,  in a nonlinear model where  a sum of solutions is not necessarily a
 solution too,  the difference between these two cases may be significant.
 Both, the function $u_{8}(S,0) $ which is a solution of equation
 (\ref{urav}) and the sum $u_3(S,0)+u_5(S,0)$ which is not equal to
 any solution of equation (\ref{urav}) are shown in Figure
 \ref{num}. We expect that if in a linear case we can use 
 linearity to compose solutions, in the nonlinear case we shall calculate
 the hedge-cost for each derivative for its own. 
Indeed, in Figure \ref{num}
 we see a strong difference between the values of the hedge-costs for the
 derivatives calculated in the linear and nonlinear cases in the
 neighbourhood of the exercise price $E.$\\

\begin{figure}
\begin{minipage}[t]{14.cm}
\includegraphics[width=14.cm]{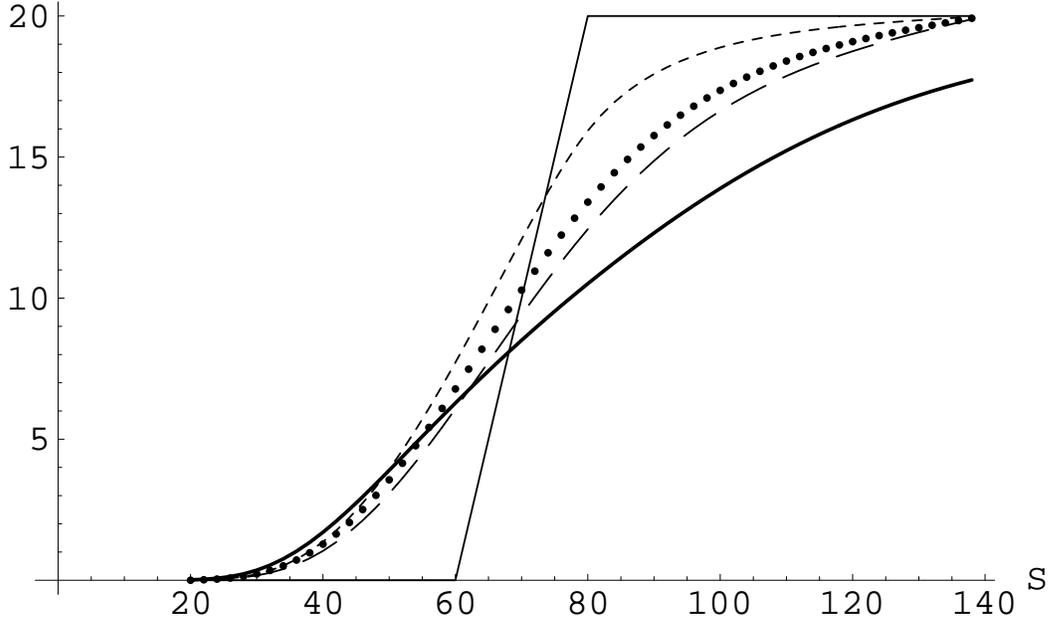}
\caption{\label{spread} 
Plot of the hedge-cost for a bull-price-spread option with the 
payoff function $u_{spread}(S,T)$, (\ref{spreadpay}), 
(thin solid line) for various values of $\rho$. Compare the
solution $u(S,0)$ of the linear Black--Scholes model (\ref{BS}) 
(thick solid line) which 
corresponds to $\rho=0$
with the numerical solutions $u_{spread}(S,0)$ to the nonlinear 
equation (\ref{urav})  
with $\rho=0.2$ (short dashed line), with $\rho=0.1$ (dots) 
and with $\rho=0.05$ (long dashed line).
 The parameters for the bull-price-spread option are $S \in (20,140]$, $t \in [0,T]$,
 $T=1.0$, $r=0.02$, 
  $\sigma=0.35$, the exercise price for the long European call is
  $E_l=60.0,$ the exercise price for the short European call is
  $E_s=80.0$.
For all numerical solutions  the same  payoff, volatility, expiry date 
and exercise prices as in linear case are chosen and 
the parameters of the grid are $ h=2,~ \tau=0.05$, $N_S=60,~ N_t=20.$}
\end{minipage}
\end{figure}

We use the completely implicit method  also for the numerical calculation of an
 hedge-cost for an option with an essential different payoff as  in
the case of a European call or a strangle. As an example we take a bull-price-spread option with the payoff 
\begin{equation}\label{spreadpay}
u_{spread}(S,T)= \max(S-E_l,0) -\max(S-E_s,0),~E_l<E_s.
\end{equation}
We used the same system of difference equations (\ref{i2}) -
(\ref{ins}) and studied the option hedge-cost for various values of
$\rho.$ The results are represented in Figure \ref{spread} and  show the strong
difference between linear and nonlinear cases of Black--Scholes
equations as well as a strong dependence of the option hedge-cost on
the feedback-effect for a large trader.

All calculations were done using the program {\bf Mathematica 5.0}.
In order to solve the system of algebraic equations (\ref{i2}) -
(\ref{ins}) we used the function {\bf FindRoot}. If we use the
boundary conditions (\ref{ten}), then on the interval $S \in [0,E]$
 $u_{bound}(S,T)=u_{K}(S,T)=0$ holds.  Would we take as the
first approximate values for the procedure {\bf FindRoot} zeros for the 
values of $u_{i,N_t-1}$ then this procedure will lead to the trivial
solution for the system of equations (\ref{i2}) - (\ref{ins}). To
avoid this problem we  take as the first approximate values for
the procedure {\bf FindRoot} some small constant $k$. We proved that
the solutions to the  system (\ref{i2}) - (\ref{ins}) do not depend on this
constant. In our calculations we used $k=0.03$ for the calculations
represented on Figures \ref{rhodep}, \ref{rhodep3}, \ref{num} and
$k=1.0$ for numerical solutions given on the Figure
\ref{spread}.\\[5pt]

\section{Conclusion}

We studied the symmetry properties of the nonlinear partial
differential equation (\ref{urav}). We found the corresponding four
dimensional Lie algebra (\ref{al4}) and the explicit representation of
the Lie group (\ref{str})--(\ref{utro}) for this equation.  The
existence of a nontrivial Lie group allowed us to obtain the invariants
(\ref{invar1})--(\ref{invar2}) which can be used as new independent
and dependent variables. Using new scaling variables we reduced the
partial differential equation to the ordinary differential equation
(\ref{uravz}). This equation possesses a solvable Lie algebra 
spanned by the infinitesimal generators
(\ref{gen}). Consequently we were able to reduce this equation to the
first order differential equation (\ref{urav_yz}). We proved
uniqueness conditions for this ordinary differential equation.  We
used the algebra (\ref{gen}) to obtain families of invariant solutions
(\ref{firu}) and (\ref{scalsol}) and proved that the uniqueness
conditions for these solutions can fail just in the point $S=0$.  The
invariant solutions have boundary conditions which approximates payoff
of strangles. We studied sensitivity parameters for these solutions
and gave graphically representations for the dependences of these
parameters on time and value of underlying asset. We used the obtained
invariant solutions to test numerical methods. We proved that the best
result can be obtained with a completely implicit method. We used
this numerical method to find numerical solutions for calls and
bull-price-spread options.  In all studied cases we have seen a strong
dependence of the option hedge-cost on the feedback-effect of the
large trader.
\section{Acknowledgments} 

 The authors are grateful to R. Frey (University of Leipzig), M. Fr{\"o}hner
    (Brandenburg University of Technology Cottbus), I.~P. Gavrilyuk (BA
    Th\"uringen), B.~N. Khoromskij (MPI Leipzig), C. Petzold (HSBC Trinkaus
    \& Burkhardt) for interesting and fruitful discussions.

The work of the second author was kindly supported by the HWP - project, grant
number 02014 of the Brandenburg,  MWFK and by  the grant of Halmstad University, Sweden.
 
\nocite{Frey:illicuidity}
\nocite{FreyPatie} 
\nocite{Frey:perfect} 
\nocite{SchonbucherWilmotta}
\nocite{Gaeta} 
\nocite{Stephani} 
\nocite{Olver} 
\nocite{Ibragimov}

\bibliographystyle{IJTAF}
\bibliography{bordagchmakova} 
\newpage

\begin{figure}[ht]
\centerline{The Greeks for  solutions (\ref{scalsol}).}
\vspace{0.5cm}
\begin{minipage}[t]{7cm}
\includegraphics[width=7cm]{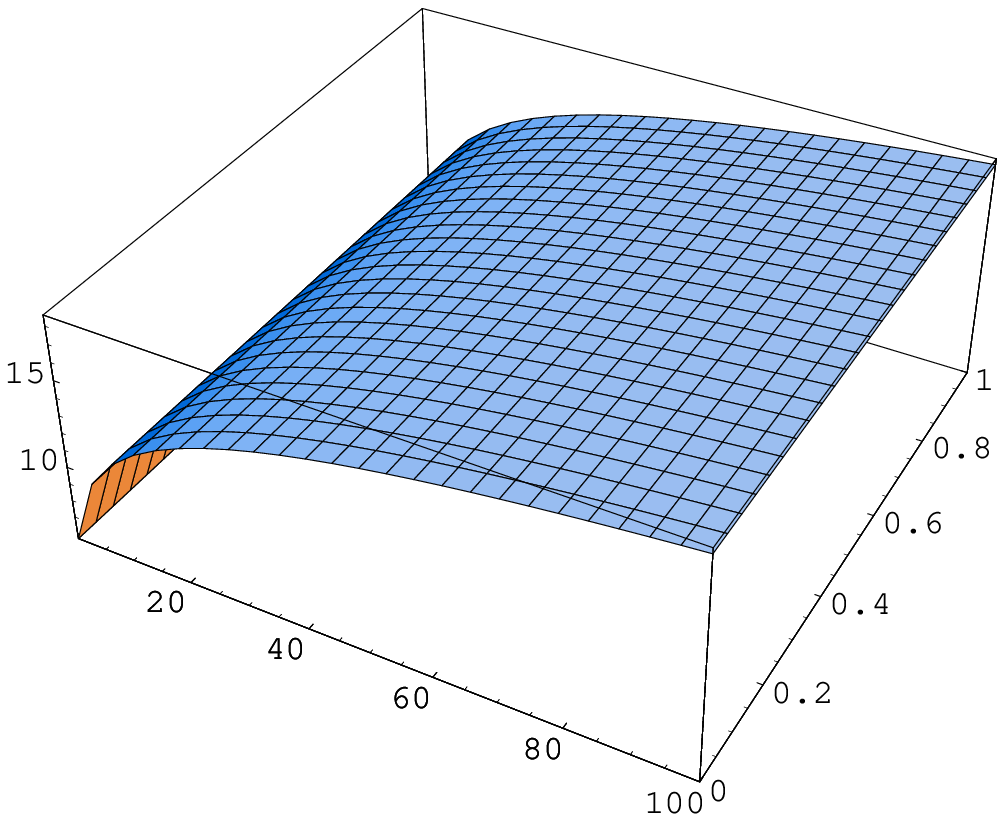}
\caption{\label{delta} Plot of
$\rho \Delta =\rho \frac{\partial u(S,t)}{\partial S}$ 
with $S \in (0,100]$, $t \in [0,1]$ and parameters 
 $\sigma=0.28,~m=8.5,~d_1=d_2=0.$}
\end{minipage} \hfill
\begin{minipage}[t]{7.cm}
\includegraphics[width=7.cm]{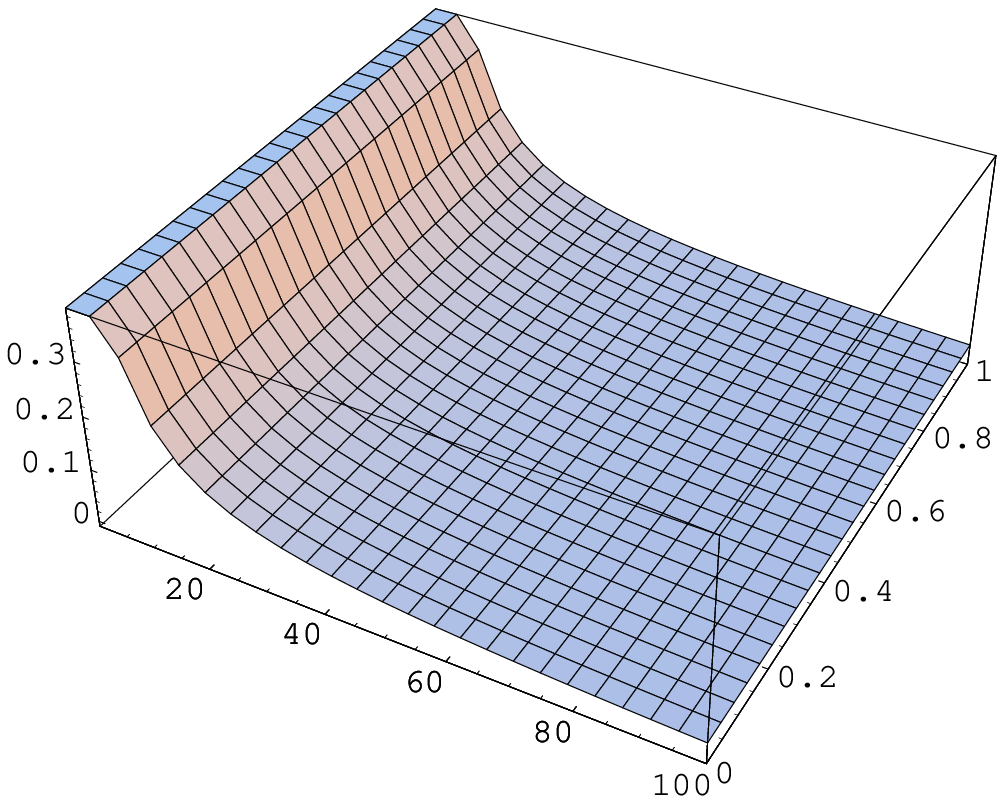}
\caption{\label{gamma} Plot of
$\rho \Gamma =\rho \frac{\partial^2 u(S,t)}{\partial S ^2}$ 
with $S \in (0,100]$, $t \in
 [0,1]$ and parameters 
$\sigma=0.35,~m=4.9$, $d_1=d_2=0.$}
\end{minipage}
\end{figure}

\end{document}